\documentclass{amsart}
\usepackage{amscd, amssymb}

\newtheorem*{remark}{Remark}

\newcommand{\refT}[1]{Theorem~\ref{T:#1}}
\newcommand{\refC}[1]{Corollary~\ref{C:#1}}
\newcommand{\refP}[1]{Proposition~\ref{P:#1}}
\newcommand{\refD}[1]{Definition~\ref{D:#1}}
\newcommand{\refL}[1]{Lemma~\ref{L:#1}}

\newcommand{\proj}{{\mathbb P}}
\newcommand{\ints}{{\mathbb Z}}
\newcommand{\cplx}{{\mathbb C}}
\newcommand{\reals}{{\mathbb R}}
\newcommand{\nats}{{\mathbb N}}
\newcommand{\wt}{\widetilde}
\newcommand{\vphi}{\varphi}
\newcommand{\n}{\underline{n}}

\newcommand{\Sf}{{\bf S}}

\newtheorem{theorem}{Theorem}[section]
\newtheorem{lemma}[theorem]{Lemma}
\newtheorem{proposition}[theorem]{Proposition}
\newtheorem{corollary}[theorem]{Corollary}
\newtheorem*{thmnonum}{Theorem}

\newtheorem{definition}[theorem]{Definition}

\begin{document}

\title{Computations of Complex Equivariant Bordism Rings}
\author{Dev Prakash Sinha}
\address{Department of Mathematics, Brown University,
Providence, RI 02906}
\email{dps@math.brown.edu}

\maketitle

\section{Introduction}

Bordism theory is fundamental in algebraic topology and its
applications.  In the early sixties Conner and Floyd introduced equivariant
bordism as a powerful tool in the study of transformation groups.  In the
late sixties, tom Dieck introduced homotopical bordism in order to refine
understanding of the localization techniques employed by
Atiyah, Segal and Singer in index theory.  Despite the many successful
computations and applications of bordism theories, equivariant bordism has    
been mysterious from a computational point of view, even for
cyclic groups of prime order $p$ (see \cite{Kos} and \cite{Kriz}).
In this paper we present the first computations of the ring structure
of the coefficients of equivariant bordism, for abelian groups.  
The key constructions are
operations on equivariant bordism.  Analogs of these operations
should play an important role in
equivariant stable homotopy more generally.  Our main techniques involve
localization and give some insight into the structure of $MU^G_*$ for
a large class of groups including $p$-groups.  We give a synopsis of our 
results now.\\

We denote by $MU^G_*$ the homotopical equivariant bordism ring, where $G$
is a compact Lie group.  It is
defined analogously to $MU_*$ as 
$\lim_V [S^{\n\oplus V}, T(\xi^G_{|V|})]^G$, where $V$ ranges over isomorphism
classes of complex representations of $G$, $S^{\n\oplus V}$ is the one-point 
compactification of the Whitney sum of $\cplx^n$ with trivial $G$ action and
$V$, and $T(\xi^G_{|V|})$ is the Thom space of the universal complex 
$G$-bundle.  In fact, we may use these Thom spaces to define an equivariant
spectrum as first done by tom Dieck \cite{tD} and hence define associated
equivariant homology and cohomology theories $MU^G_*(-)$ and $MU^*_G(-)$.
We will carefully make these constructions in section \ref{S:MUG}.

Euler classes play fundamental roles in our work.  The Euler classes which
are most important for us are those associated to a complex
representation of $G$, considered as a $G$-bundle over a point.  
More explicitly,
the Euler class associated to $V$ is a class $e_V \in MU_G^m(pt.)$, where $m$
is the dimension of $V$ over the reals, represented 
by the composite $S^0 \hookrightarrow S^V \to T(\xi^G_{|V|})$,
where the second map is ``inclusion of a fiber''.  Euler classes
multiply by the rule $e_V \cdot e_W = e_{V \oplus W}$.
In homological grading $e_V \in MU^G_{-m}$, so it
cannot be in the image of a geometric bordism class under
the Pontrijagin-Thom map if it is non-trivial.  
If $V^G = \{0\}$ then $e_V$ is non-zero, reflecting
the fact that $V$ has no non-zero equivariant sections.
Therefore, the homotopy groups of $MU^G$ are not bounded below, 
a feature which already distinguishes it from its ordinary counterpart.

More familiar classes in $MU^G_*$ are those in the image of
classes in geometric bordism under the Pontrijagin-Thom map.  
Given a stably complex  $G$-manifold $M$, let $[M]$ denote the corresponding
class in $MU^G_*$.  Complex projective spaces give a rich collection
of examples of $G$-manifolds.
Given a complex representation $W$ of $G$ let $\proj(W)$ denote the space of
complex one-dimensional subspaces of $W$ with inherited $G$-action.

The starting point in our work is that after inverting Euler classes, 
$MU^G_*$ becomes computable by non-equivariant means.  
That we rely heavily on localization is not surprising
because localization techniques have pervaded equivariant topology.  For
any compact Lie group $G$ let
$R_0$ denote the sub-algebra of $MU^G_*$ generated by the $e_V$ and 
$[\proj(\n \oplus V)]$ as $V$ ranges over non-trivial irreducible 
representations.  Let $S$ be the multiplicative set in $R_0$  
of non-trivial Euler classes.  By abuse, denote the same multiplicative 
set in $MU^G_*$ by $S$.  Then the key first result, which we emphasize
is true for a large class of groups including $p$-groups, is the following.

\begin{theorem}\label{T:locisom}
Let $G$ be a group such that any proper subgroup is contained in a proper
normal subgroup.  The inclusion of $R_0$ into $MU^G_*$ becomes
an isomorphism after inverting $S$.
\end{theorem}

In other words, we may multiply any class in $MU^G_*$ by some Euler class
to get a class in $R_0$ modulo the kernel of the localization map $S$.
We are lead to study divisibility by Euler classes as well as
the kernel of this localization map.  We can do so successfully
in the case when the group in question is a torus.

Let $T$ be a torus, and let $V$ be a non-trivial irreducible representation
of $T$.  Let $K(V)$ denote the subgroup of $T$ which acts trivially on $V$.
There is a restriction homorphism (of algebras) 
$res^T_H \colon MU^T_* \to MU^H_*$ for any subgroup $H$.
The restriction of $e_V$ to $MU^{K(V)}_*$ is zero, as can be seen using an
explicit homotopy.  Remarkably, we have the following.

\begin{theorem}\label{T:seq}
The sequence
$$ 0 \to MU^T_* \overset{\cdot e_V}{\to} MU^T_* \overset{res^T_{K(V)}}{\to}
MU^{K(V)}_* \to 0 $$
is exact.
\end{theorem}

Note that the surjectivity of the restriction map is false for geometric 
bordism.  One cannot for example extend the non-trivial action of $\ints/2$
on two points to an $S^1$-action.

Using this exact sequence, we define operations which are essentially 
division by Euler classes.  To define these operations we need to split the
restriction maps.  The restriction map to the trivial group
is called the augmentation map $\alpha \colon MU^G_* \to MU_*$.  There 
is a canonical splitting of this map as rings which defines an $MU_*$-algebra
structure on $MU^G_*$.  All of the maps we have defined so far are in fact
maps of $MU_*$-modules.  The restriction maps to other sub-groups are not
canonically split, but we do know the following from \cite{CCM}.

\begin{theorem}[Comeza\~na]\label{T:comez}
Let $G$ be abelian.  Then $MU^G_*$ is a free $MU_*$-module concentrated
in even degrees.
\end{theorem}

Hence we may fix a splitting $s_V$ as $MU_*$-modules
of the restriction map $res^T_{K(V)}$.  Unless
$K(V)$ is the trivial group, this splitting is non-canonical and is not
a ring homomorphism.

\begin{definition}
Let $T$ and $V$ be as above.  Define the $MU_*$-linear operation
$\Gamma_V$ as follows.
Let $x \in MU^T_*$.  Then $\Gamma_V(x)$ is the unique class in $MU^T_*$ 
which satisfies $$e_V \cdot \Gamma_V(x) = x - s_V(res^T_{K(V)} x).$$
\end{definition}

For convenience, let $\beta_V$ denote $s_V \circ res^{T}_{K(V)}$.
We are now ready to state our main theorem. 

\begin{theorem}\label{T:main}
For any choice of splittings $s_V$, $MU^T_*$ is generated as an $MU_*$-algebra 
over the operations $\Gamma_V$ 
by the classes $e_V$ and $[\proj(\n \oplus V)]$, where $n \in \nats$ and
in all instances $V$ ranges over non-trivial
irreducible  complex representations of $T$.
Relations are as follows:
\begin{itemize}
\item $e_V \Gamma_V(x) = x - \beta_V(x)$
\item $\Gamma_V(\beta_V(x)) = 0$ 
\item $\Gamma_V(e_V) = 1$ 
\item $\Gamma_V(xy) = \Gamma_V(x) y + \beta_V(x) \Gamma_V(y) + \Gamma_V(\beta_V(x) \beta_V(y))$
\item $ \Gamma_V \Gamma_W x  = \Gamma_W \Gamma_V x   + 
\Gamma_W \Gamma_V \beta_W(x)  - \Gamma_W \Gamma_V (e_W) \beta_V(\Gamma_W x) \\
\;\;\;\;\;\;\;\;\;\;\;\;\;\;\;\;\; - \Gamma_W \Gamma_V (\beta_V(e_W) \beta_V(\Gamma_W x)),$

\end{itemize}
where $V$ and $W$  range over non-trivial irrreducible representations of
$T$ and $x$ and $y$ are any classes in $MU^G_*$.
\end{theorem}

We may recover the structure of $MU^G_*$ for any abelian group $G$ by
realizing $G$ as the kernel of an irreducible representation of some 
torus and using the exact sequence of \refT{seq}.

%The sub-algebra $R_0$ of $MU^T_*$ 
%is a maximal polynomial sub-algebra.  Hence the finer structure of $MU^T_*$,
%which has been a notoriously difficult algebra to describe, is embodied
%by the operations $\Gamma_V$.  We call these Conner-Floyd
%operations because in the special case of geometric classes
%where $T$ is $S^1$ and $\rho$ is the standard representation
%there is a geometric model for $\Gamma_\rho([M])$ due to Conner and Floyd.

We give both algebraic and geometric applications of our main computation.
For $G = S^1$ and $\rho$ its standard representation, 
we present a geometric model of $\Gamma_\rho([M])$.  This geometric model 
allows us to compute the 
completion map $MU^G_* \to (MU^G_*)_{\hat{I}}$, where $I$ is the kernel 
of the augmentation map from $MU^G_*$ to $MU_*$.  The completion theorem
of L\"offler, as proved by Comeza\~na and May, states that for $G$ abelian,
$(MU^G_*)_{\hat{I}} \cong MU^*(B_G)$, where $B_G$ is the classifying space of 
$G$.  So this completion map gives a connection between equivariant
bordism and any equivariant theory which is defined using a Borel
construction $E_G \times_G -$.
We also give more classically-styled 
applications to the understanding of group actions on manifolds.
For example, a current
topic of great interest in equivariant cohomology is the investigation of 
$G$-manifolds with isolated fixed points, essentially extending Smith theory.
We prove the following.

\begin{theorem}\label{T:appl}
Let $M$ be a stably-complex four dimensional $S^1$-manifold with three 
isolated fixed points.   Then $M$ is equivariantly cobordant to 
$\proj(\underline{1}\oplus V \oplus W)$ for some distinct non-trivial
irreducible representations $V$ and $W$ of $S^1$.
\end{theorem}

%The structure of this paper is as follows.  In the next two sections we 
%review basic definitions.  In section~4 we establish the well-known connection
%between taking fixed sets and localization in equivariant homotopy theory,
%in our setting of equivariant bordism.  Section~5 is brief, introducing
%the algebraic language which has proven to be most useful in equivariant
%bordism theory.
%We prove our main theorems in section~6. 
%We give geometric applications in section~7, answering questions posed by
%Bott and Milgram.

The author thanks his thesis advisor, Gunnar Carlsson, for pointing him to 
this problem and for innumerable helpful comments.  As this project has spanned
a few years, the author has many people to thank for conversations which
have been helpful including Botvinnik, Goodwillie, Klein, Milgram, Sadofsky, 
Scannell, Stevens and Weiss.  He also thanks Haynes Miller for a close
reading of an earlier version of this paper.  Thanks also go to 
Greenlees, Kriz and May for sharing preprints of their work.  

\section{Preliminaries}

Until otherwise noted, the group $G$ is a compact Lie group. 

All $G$ actions are assumed to be continuous, and $G$-actions on
manifolds are assumed to be smooth.  For any $G$-space X, we let $X^G$
denote the subspace of $X$ fixed under the action of $G$.  The space 
of maps between two $G$-spaces, which we denote ${\text{Maps}}(X, Y)$
has a $G$-action by conjugation.
We denote its subspace of $G$-fixed maps by ${\text{Maps}}^G(X, Y)$.
We will often work with based spaces, in which case we assume that the
basepoints are fixed by $G$.  
Throughout,
$E_G$ will be a contractible space on which $G$ is acting freely.  And
$B_G$, the classifying space of $G$, is the quotient of $E_G$ by the
action of $G$.

We will always let $V$ and $W$ be finite-dimensional complex
representations of $G$.  
Our $G$-vector bundles will always have paracompact base
spaces, so we may define a $G$-invariant inner product on the fibers.
The constructions we make using such an inner product will be independent
of choice of inner product up to homotopy.  
We will use the same notation for a $G$-bundle over a point
as for the corresponding representation.  We let $|V|$ denote the
dimension of $V$ as a complex vector space.  The sphere $S^V$ is the
one-point compactification of $V$, based at $0$ if a base point is
needed.  And the sphere $S(V)$ is the unit sphere in $V$ with
inherited $G$-action.  For a $G$-vector bundle $E$, let $T(E)$ denote
its Thom space, which is the cofiber of the unit sphere bundle of $E$
included in the unit disk bundle of $E$.  Thus for $V$ a
representation $T(V) = S^V$.

Let $R^+(G)$ denote the monoid (under direct sum) of isomorphism
classes of complex representations of $G$, and let $R(G)$ denote the
associated Grothendieck ring (where multiplication is given by tensor
product).  We let $Irr(G)$ denote the set of isomorphism classes of 
irreducible complex representations of $G$, and let $Irr^*(G)$ be the
subset of non-trivial irreducible representations.  
If $W = \sum a_i V_i \in R(G)$ where
$V_i$ are distinct irreducible representations, let $\nu_V(W)$ for an 
irreducible $V$ be $a_j$
if $V$ is isomorphic to $V_j$ or zero if $V$ is not isomorphic to any
of the $V_i$.  Recall from
the introduction that $\rho$ is the standard representation of $S^1$.
We will by abuse use $\rho$ to denote the standard representation
restricted to any subgroup of $S^1$.  We use $\underline{n}$ or $\cplx^n$
to denote the trivial $n$-dimensional complex representation of a
group.  We will sometimes think of representations as group homomorphisms,
and talk of their kernels, images, and so forth.

We rely on techniques from equivariant stable homotopy theory.
Let $\Omega^W(X)$ denote the space of based maps from $S^W$ to $X$.  
By fixing a representation $\mathcal{U}$ with inner product, 
of which a countably infinite direct sum of any
representation of $G$ appears as a summand, we define a $G$-spectrum $X$ to 
be a family of spaces $X_V$ indexed on subspaces of $\mathcal{U}$
equipped with $G$-homeomorphisms $X_V \to \Omega^{W \ominus V} X_W$
for all $V \subseteq W$.  The basic passage to ordinary stable homotopy
theory is by taking the fixed-points spectrum.  Consider only subspaces
$V \subset \mathcal{U}^G$.  Then we may define the fixed-points
spectrum $X^G$ using the 
family of spaces $(X_V)^G$, where the bonding maps are restrictions to
fixed sets of the given bonding maps.

\section{Basic Properties of $MU^G$}\label{S:MUG}

There are two basic definitions of bordism, geometric and homotopy
theoretic.  Equivariantly, these two theories are not equivalent, and we
will comment on this difference later in this section.

Our main concern is the homotopy theoretic version of complex
equivariant bordism, as first defined by tom Dieck \cite{tD}.  Fix
$\mathcal U$, a complex representation of which a countably
infinite direct sum of any
representation of $G$ appears as a summand.  
If there is ambiguity
possible we specify the group by writing ${\mathcal U}(G)$.  
Let $BU^G(n)$ be the Grassmanian of complex
$n$-dimensional linear subspaces of $\mathcal U$.  
Let $\xi^G_n$ denote the tautological complex $n$-plane bundle
over $BU^G(n)$.  As in the non-equivariant setting, the bundle $\xi^G_n$
over $BU^G(n)$ serves as a model for the universal complex $n$-plane
bundle.  If $V$ is a complex representation, set $\xi^G_V = \xi^G_|V|$.

\begin{definition}\label{D:MUG} 
We let $TU^G$ be the pre-spectrum, indexed on all complex
subrepresentations of $\mathcal U$, defined by taking the $V$th entry to be
$T(\xi^G_{V})$ (it suffices to define entries of a prespectrum only
for complex representations).  Define the bonding maps by noting that
for $V \subseteq W$ in $\mathcal U$, letting $V^\bot$ denote the
complement of $V$ in $W$, we have
$$S^{V^\bot} \wedge T(\xi^G_{V}) \cong T({V^\bot} \times \xi^G_{V}).$$
Then use the classifying map
$${V^\bot} \times \xi^G_{V} \to \xi^G_{W}$$ to define the
corresponding map of Thom spaces.  Pass to a spectrum in the usual
way, so that the $V$th de-looping is given by
$$\lim_{W \supseteq V} \Omega^{V^\bot}(T(\xi^G_{W})),$$ to obtain the
homotopical equivariant bordism spectrum $MU^G$.
\end{definition}

From this spectrum indexed by subspaces of $\mathcal{U}$ we may pass to
an $RO(G)$-graded homology theory $MU^G_{\star}(-)$.  We will be concerned 
with the coefficient ring in integer gradings, which we denote $MU^G_*$.
But for some arguments, we will need groups graded by complex representations
of $G$, giving rise to the need for the following proposition.

\begin{proposition}\label{P:period}
Let $V$ be a complex representation of $G$.
The group $MU^G_{V}(X)$ is naturally isomorphic to $MU^G_{2|V|}(X)$.
\end{proposition}

We prove this proposition after defining the needed multiplicative 
structure on $MU^G$.  The classifying map of the Whitney sum
$$\xi^G_{V} \times \xi^G_{W} \to \xi^G_{V \oplus W}$$ gives rise to a
map
$$T(\xi^G_{V}) \wedge T(\xi^G_{W}) \to T(\xi^G_{V \oplus W}),$$ which
defines a multiplication on $MU^G$.  The unit element is represented
by the maps $S^V \to T(\xi^G_{V})$ induced by passing to Thom spaces
the classifying map of $V$ viewed as a $G$-bundle over a point.  Thus
in the usual way the coefficients $MU^G_{\star}$ form a ring and
$MU^G_{\star}(X)$ is a module over $MU^G_{\star}$.  

\begin{definition}\label{D:thomclass}
Let $V \subset \mathcal{U}$ be of dimension $n$.  Then the classifying map $V
\to \xi^G_n$ induces a map of Thom spaces $S^V \to T(\xi^G_{n})$, which
represents an element $t_V \in MU^G_{V-2n}$ known as a {\em Thom
class}.
\end{definition}

\begin{proof}[Proof of \refP{period}]
We show that the Thom class $t_V$ is invertible.  The isomorphism between
$MU^G_{V}(X)$ and $MU^G_{2|V|}(X)$ is then given by multiplication by this
Thom class.

The class in $MU^G_{2n - V}$ represented by the map $S^{2n} \to T(\xi^G_{V})$
induced by the classifying map $\cplx^n \to \xi^G_{V}$ is the 
multiplicative inverse of $t_V$.  The product of this class with $t_V$ is
homotopic to the unit map $S^{V \oplus \cplx^n} \to 
T(\xi^G_{V \oplus \cplx^n})$.
\end{proof}

The most pleasant way to produce classes in $MU^G_*$ is from equivariant
stably almost complex manifolds.
Recall that there is an real analog of $BU^G(n)$, which we call
$BO^G(n)$, and which is the classifying space for all $G$-vector bundles.

\begin{definition}
A tangentially complex $G$-manifold is a pair $(M, \tau)$ where $M$ is a
smooth $G$-manifold and $\tau$ is a lift to $BU^G(n)$ of the map to
$BO^G(2n)$ which classifies $TM \times \reals^k$ for some $k$.
\end{definition}

We can define bordism equivalence in the usual way to get a geometric
version of equivariant bordism.

\begin{definition}
Let $\Omega^{U, G}_*$ denote the ring of tangentially complex $G$-manifolds
up to bordism equivalence.
\end{definition}

Classes in geometric bordism give rise to classes in homotopical bordism
through the Pontrjagin-Thom construction.

\begin{definition}
Define a map $PT\colon \Omega^{U,G}_* \to MU^G_*$ as follows.  Choose a
representative $M$ of a bordism class.  Embed $M$ in some sphere $S^V$,
avoiding the basepoint and
so that the normal bundle $\nu$ has a complex structure.  Identify the
normal bundle with a tubular neighborhood of $M$ in $S^V$.
Define $PT([M])$ as the composite
$$ S^V \overset{c}{\to} T(\nu) \overset{T(f)}{\to} T(\xi_{|\nu|}), $$
where $c$ is the collapse map which is the identity on $\nu$ and sends
everything outside $\nu$ to the basepoint in $T(\nu)$, and $T(f)$ is the
map on Thom spaces given rise to by the classifying map $\nu \to \xi_{|\nu|}$.
\end{definition}

The proof of the following theorem translates almost word-for-word from
Thom's original proof.

\begin{theorem}\label{T:PT}
The map $PT$ is a well-defined graded ring homomorphism.
\end{theorem}

The Pontrjagin-Thom homomorphism is not an isomorphism equivariantly 
as it is in the ordinary setting.  A theorem of Comeza\~na
states that $PT$ is split injective for abelian groups.  The following
classes illustrate the failure of the Pontrijagin-Thom map to be an 
isomorphism.

\begin{definition}\label{D:euler}
Compose the map $S^V \to T(\xi^G_{n})$, in \refD{thomclass} of the Thom
class with the evident inclusion $S^0 \to S^V$ to get an element $e_V
\in MU^G_{-2n}$ which is called the {\em Euler class} associated to
$V$.
\end{definition}

We will see that Euler classes $e_V$ associated to representations $V$ such
that $V^G = \{ 0 \}$ are non-trivial.  Thus
$MU^G_*$ is not connective, a feature which already distinguishes it
from $\Omega^{U,G}_*$ as well as $MU_*$.  The key difference
between the equivariant and ordinary settings is the lack of transversality
equivariantly.  For example, if $V^G = \{ 0 \}$ the inclusion of
$S^0$ into $S^V$ cannot be deformed equivariantly to be transverse regular
to $0 \in S^V$.

Finally, we introduce maps relating bordism rings for different groups.  
Recall that ordinary homotopical bordism $MU$
can be defined using Thom spaces as in our definition of $MU^G$ but
without any group action present.

\begin{definition}
Define the augmentation map $\alpha\colon MU^G \to MU$ by forgetting
the $G$-action on $MU^G$.  When $G$ is abelian and $H$ is a subgroup of $G$
define ${res}^G_H$ to be the map from $MU^G_* \to MU^H_*$ by restricting
the $G$-action to an $H$-action.  
\end{definition}

We need to have $G$ abelian for the map ${res}^G_H$ to be so defined.
In the abelian setting, any complex
representation of $H$ extends to a complex
representation of $G$, so that when its $G$-action is
restricted to an $H$-action the
Thom space $T(\xi^G_n)$ coincides with $T(\xi^H_n)$.

\begin{definition}\label{D:incl}
Define the inclusion map $\iota\colon MU \to MU^G$ by composing a map
$S^n \to T(\xi_n)$ with the inclusion $T(\xi_n) \to T(\xi^G_n)$.
%More generally, 
%we may define an inclusion from $MU^G$ to any $MU^{G \times H}_*$
%by imposing a trivial $H$-action on a $G$-map from a sphere to 
%Thom space and including $T(\xi^G_n)$ into the $T(\xi^{G\times H}_n)$ 
%by including ${\mathcal U}(G)$ as the $H$-fixed set of
%${\mathcal U}(G \times H)$ and passing to Thom spaces.
\end{definition}

On coefficients, $\iota$ defines an $MU_*$-algebra structure on
$MU^G_*$.  The kernel of $\alpha$ on coefficients is called the {\em
augmentation ideal}.  For example, the Euler class $e_V$ is in the
augmentation ideal as the map $S^0 \to S^V$ in its definition is
null-homotopic when forgetting the $G$-action.
%\begin{definition}
%Let $\bar{x}$ denote the image of $x \in MU^G_*$ under $\iota \circ
%\alpha$.
%\end{definition}
%
%Then the augmentation ideal contains, and is clearly generated by,
%elements of the form $x - \bar{x}$.
On the other hand, $\iota$ is injective, which follows from the
following proposition which is proved for example in \cite{CCM}.

\begin{proposition}\label{P:split1}
The composite $\alpha \circ \iota\colon MU \to MU$ is homotopic to the
identity map.
\end{proposition}

\section{The Connection Between Taking Fixed Sets 
                      and Localization}\label{S:loc1}

The connection between localization, in the commutative algebraic
sense, and ``taking fixed sets'' has been a fruitful theme in
equivariant topology.  We develop this connection in the setting
of bordism in this section.

The main goal of this section is to prove \refT{locisom}, which we
restate here for convenience.  Let
$R_0$ denote the sub-algebra of $MU^G_*$ generated by the classes $e_V$ and 
$[\proj(\n \oplus V)]$ as $V$ ranges over non-trivial irreducible 
representations.  
Let $S$ be the multiplicative set in $R_0$  
of non-trivial Euler classes.

\begin{definition}
A compact Lie group $G$ is {\em nice} if every proper subgroup is
contained in a proper normal subgroup.
\end{definition}

For example, abelian groups and $p$-groups are nice.

\begin{thmnonum}[Restatement of \refT{locisom}]
Let $G$ be a nice group.  The inclusion of $R_0$ into $MU^G_*$ becomes
an isomorphism after inverting $S$.
\end{thmnonum}

We prove this theorem by first explicitly computing $S^{-1}MU^G_*$ and
then computing the images of generators of $R_0$ in $S^{-1}MU^G_*$.
We start with the following lemma, 
which provides translation between localization and topology.  For
any commutative ring $R$ and element $e \in R$ let $R[\frac{1}{e}]$
denote the localization of $R$ obtained by inverting $e$.

\begin{lemma}\label{L:inveul}
As rings, 
$\widetilde{MU^G_*}(S^{\oplus_\infty V}) \cong MU^G_*[\frac{1}{e_V}].$
\end{lemma}

\begin{proof}
The left-hand side $\widetilde{MU^G_*}(S^{\oplus_\infty V})$ is a ring
because $S^{\oplus_\infty V}$ is an $H$-space via the equivalence
$$S^{\oplus_\infty V} \wedge S^{\oplus_\infty V} \cong S^{\oplus_\infty V}.$$
To compute the left-hand side, 
apply $\wt{MU^G_*}$ to the identification $S^{\oplus_\infty V} =
\varinjlim S^{\oplus_n V}$.  After applying the suspension
isomorphisms $\wt{MU^G_*}(S^{\oplus_k V}) \cong
\wt{MU^G_{*+|V|}}(S^{\oplus_{k+1}V})$, the maps in the resulting
directed system are multiplication by the $e_V$.
\end{proof}

We will see that after inverting Euler classes, 
equivariant bordism is computable.
If $G$ is a nice group and $\{ W_i \}$ are the non-trivial irreducible
representations of $G$ then $Z = S^{\oplus_i(\oplus_\infty W_i)}$ has fixed
sets $Z^G = S^0$ while $Z^H$ is contractible for any $H \subset G$.
Hence, our next lemma, taken with \refL{inveul}, 
establishes the strong link between
localization and taking fixed sets.

\begin{lemma}\label{L:l1}
Let $X$ be a finite $G$-complex and let $Z$ be a $G$-space such that
$Z^G \simeq S^0$ and $Z^H$ is contractible for any proper subgroup of
$G$.  Then the restriction map
$${\text{Maps}}^G(X, Y \wedge Z) \to {\text{Maps}}(X^G, Y^G)$$ 
is a homotopy equivalence.
\end{lemma}

\begin{proof}
The restriction map is a fibration whose fiber at a given point
is the space of $G$-maps which are specified on $X^G$.
Using the skeletal filtration of $X$, we can then filter this mapping 
space by spaces 
$${\text{Maps}}^G(D^k \times {G/H}, Y \wedge Z),$$
such that the maps are specified on the boundary of $D^k \times {G/H}$, and
where $H$ is a proper subgroup of $G$.  A standard change-of-groups argument
yields that this mapping space is homeomorphic to 
${\text{Maps}}(D^k , (Y\wedge Z)^H)$, again with the map specified on the 
boundary.   But $(Y \wedge Z)^H$ is contractible, and thus so are
these mapping spaces.  Thus, the fiber of the restriction map is 
contractible.
\end{proof}

We now translate this lemma to the stable realm.  For
simplicity, let us suppose that our $G$-spectra are indexed over the
real representation ring.  We can do so by choosing specific representatives
of isomorphism classes of representations.
Let $K_n \subset K_{n+1}$ denote a
sequence of representations which eventually contain all irreducible
representations infinitely often and such that ${K_n}^\perp \subset
K_{n+1}$ contains precisely one copy of the trivial representation.
If $G$ is finite, we can let $K_n$ be the direct sum of $n$ copies of
the regular representation.

\begin{definition}
Let $X$ be a $G$-prespectrum.  We define the geometric fixed sets
spectrum $\Phi^G X$ by passing from a prespectrum $\phi^G X$ defined
as follows.  We let the entry $\{\phi^G X\}_n$ be $(X_{K_n})^G$, the
$G$-fixed set of the $K_n$-entry of $X$.  The bonding maps are
composites
$$
(X_{K_n})^G {\longrightarrow} (\Omega^{{K_n}^\perp}X_{K_{n+1}})^G 
\longrightarrow
\Omega^{({K_n}^\perp)^G} (X_{K_{n+1}})^G = \Omega (X_{K_{n+1}})^G,
$$
where the first map is a restriction of a bonding map of $X$, and the 
second map is restriction to fixed sets of the loop space.
\end{definition}

While the prespectrum $\phi^G(X)$ depends on the choice of filtration
$K_*$, the spectrum $\Phi^G X$ is independent of this choice.

\begin{lemma}\label{L:PhiG}
Let $Z$ be as in \refL{l1}.  Then for any $G$-prespectrum $X$, the
prespectra $(X \wedge Z)^G$ and $\Phi^G X$ are homotopy equivalent.
\end{lemma}

\begin{proof}
From the definition of $(X \wedge Z)^G$, consider
$$(\Omega^W(X_{W \oplus V} \wedge Z))^G.$$ Applying \refL{l1}, 
the restriction from this mapping space to
$\Omega^{W^G} (X_{W \oplus V})^G$ is a homotopy equivalence.  Choosing
$V = K_n$, we see that $\Omega^{W^G} (X_{W \oplus K_n})^G$ is an entry
of $\phi^G X$.  The bonding maps clearly commute with these
restriction to fixed sets maps, so we have an equivalence of spectra.
\end{proof}

Note that any $Z$ as in \refL{l1} is an (equivariant)
H-space as $Z \wedge Z \simeq Z$. 
Hence if $X$ is a ring spectrum so is $(X \wedge Z)^G$.
Taking \refL{inveul} and \refL{PhiG} together, we have the following.

\begin{proposition}\label{P:blah}
Let $G$ be a nice group and let $S$ be the multiplicative set of 
non-trivial Euler classes in $MU^G_*$.  Then as rings
$S^{-1}MU^G_* \cong (\Phi^G MU^G)_*$.  
\end{proposition}

To compute $(\Phi^G MU^G)_*$, we can use the geometry of Thom spaces.  
Because smashing a weak equivalence of prespectra with a complex yields
another weak equivalence, we have $Z \wedge MU^G \simeq Z \wedge TU^G$ as 
prespectra, where  $TU^G$ denotes the equivariant Thom
prespectrum and $Z$ is as in \refL{l1}.  Hence,
$$
   \Phi^G MU^G \simeq Z \wedge MU^G \simeq Z \wedge TU^G \simeq \Phi^G
        TU^G,
$$
As required by the definition of $\Phi^G$, we
proceed with analysis of fixed-sets of Thom spaces.  

We first need the following basic fact about equivariant vector bundles.

\begin{proposition}\label{P:decomp}
And let $E$ be a $G$-vector bundle over a base space with trivial $G$-action
$X$.  Then $E$ decomposes as a direct sum
$$ E \cong \bigoplus_{V \in Irr(G)} E_{V},$$ 
where $E_V \cong \wt{E} \otimes V$ for some vector bundle $\wt{E}$.
\end{proposition}

The following result is due to tom Dieck \cite{tD} in the case of the group
$G = \ints/p$.

\begin{lemma}\label{L:thomfs}
For any compact Lie group $G$, 
the $G$-fixed set of the Thom space of $\xi^G_n$ is homotopy equivalent to
$$
\bigvee_{W \in {R^+(G)}_n} T({\xi}_{|W^G|}) \wedge 
\left(\prod_{V \in Irr^*(G)} 
BU(\nu_V(W) ) \right)_+,
$$
where we define $R^+(G)_n$ as the subset of dimension $n$ representations 
in $R^+(G)$ and we recall that $\nu_V(W)$ is the greatest number $m$ such
that $\oplus_m V$ appears as a summand of $W$.
\end{lemma}

\begin{proof}
The universality of $\xi^G_n$ implies that $(BU^G(n))^G$ is a
classifying space for $n$-dimensional complex $G$-vector bundles over
trivial $G$-spaces.  Using \refP{decomp} we see that this classifying
space is weakly equivalent to
$$\coprod_{W \in R^+(G)} \left( \prod_{V \in Irr(G)} 
BU(\nu_V(W) ) \right).$$ 
Over each component of
this union, the universal bundle decomposes as $\xi_1 \times \xi_2$,
where $\xi_1$ is the universal vector bundle over the factor of 
$\prod BU(n)$ corresponding to the trivial representation.
The fixed set ${\xi_1}^G$ is all of $\xi_1$ 
while the fixed set ${\xi_2}^G$ is the zero
section.  The result now follows by passing to Thom spaces.
\end{proof}

For convenience, we define the following spectrum.  

\begin{definition}
Let $$\Sf[Irr^*(G)] = \bigvee_{W \in \ints[Irr^*(G)] 
\subset R(G)} S^{2(|W|)}.$$
Define a ring spectrum structure on $\Sf[Irr^*(G)]$ by sending the
$V$ summand smashed with the $W$ summand to the $V+W$ summand.
\end{definition}

\begin{theorem}\label{T:phiGspec}
For any compact Lie group $G$, 
$$\Phi^G MU^G \simeq \Sf[Irr^*(G)] \wedge MU 
\wedge (\prod_{V \in Irr^*(G)} BU)_+.$$
\end{theorem}

After \refL{thomfs}, the proof of this theorem is straightforward, passing
from the prespectrum $\phi^G TU^G$ to the spectrum $\Sf[Irr^*(G)] \wedge MU 
\wedge (\prod_{V \in Irr^*(G)} BU)_+.$

For a non-trivial irreducible representation $V$, let $f_V$ be the map
from $\cplx \proj^k$ mapping to the $V$th wedge summand of
$\prod_{V \in Irr^*(G)} BU$ by the
canonical inclusion to $BU(1) \subset BU$ on the $V$th factor and by the
trivial map on the other factors. 
Define $Y_{i, V}$ to be 
the class in the subgoup $MU_{2(i-1)} (\prod_{V \in Irr^*(G)} BU)_+)$ of
$(\Phi^G MU^G)_{2i}$ represented by $\cplx \proj^{i-1}$

We may now complete the central computation of this section.

\begin{theorem}\label{T:floc}
The ring $(\Phi^G MU^G)_*$ is a Laurent algebra tensored with a 
polynomial algebra as follows:
$$(\Phi^G MU^G)_*
\cong MU_* \left[ e_V, e_V^{-1}, Y_{i, V} \right].$$   Here $V$ ranges over
irreducible representations of $G$, $i$ ranges over the positive
integers, where as indicated by notation $e_V$ is the image of 
the Euler class $e_V \in MU^G_*$ under the canonical map to the 
localization and where $Y_{i, V}$ are as above.
\end{theorem}

\begin{proof}
This theorem is simple computation after \refT{phiGspec}. 
We use the computation $MU_*(BU) \cong MU_*[Y_i]$ as rings, where
$Y_i$ is represented by $\cplx\proj^i$ mapping to $BU$ via its inclusion
into $BU(1)$,
which is standard as in \cite{Ad}.  Because $MU_*(BU)$ is a free
$MU_*$-module, it follows from the K\"unneth theorem
that $MU_*(\prod_{Irr^*(G)} BU)$ is a polynomial algebra as well.
To finish the computation, we note that the Euler class $e_V$ maps to the
class in  of $\pi_{-|V|} \Sf[Irr^*(G)]$ which is the generator on the
$V$th summand.  
\end{proof}

From \refP{blah} and the above theorem we have the following.

\begin{corollary}\label{C:comploc}
$S^{-1}MU^G_* \cong MU_*[e_V, e_V^{-1}, Y_{i,V}].$
\end{corollary}

We have shown the intimate relation between localization and taking 
fixed sets for homotopical equivariant bordism.
We will also need the following geometric point of view, which dates
back to Conner and Floyd.  

\begin{proposition}\label{P:fixed}
Let $M$ be a tangentially complex $G$-manifold.
The normal bundle $\nu$ of $M^G$ in $M$ is a complex vector bundle.
\end{proposition}

\begin{proof}
Let $\eta$ be a complex $G$-bundle over $M$ whose underlying real bundle
is $TM \times \reals^k$, as given by the tangential unitary $G$-structure of
$M$.  Then by \refP{decomp}, $\eta |_{M^G}$ decomposes
as a complex $G$-bundle
$$\eta|_{M^G} \cong \eta_1 \oplus
                \bigoplus_{\rho \in Irr^0(G)} \eta_{\rho},$$
where $\eta_1$ has trivial $G$-action.  But we can identify $\eta_1$ as having
underlying real bundle equal to $TM^G \times \reals^k$.  So the normal
bundle $\nu$ underlies $\bigoplus_{\rho \in Irr^0(G)} \eta_{\rho}$, which
gives $\nu$ the desired complex structure.
\end{proof}

As Comeza\~na points out, 
this proposition would not be true if in the 
defintion of complex $G$-manifold we chose a complex structure on either
the stable normal bundle or on $TM \times V$ for an arbitrary $V$ as
opposed to $\reals^k$.  In these cases we could only guarantee that normal
bundles to fixed sets would be stably complex.

\begin{definition}\label{D:vphi}
Let 
$$F_* = \bigoplus_{W \in R^+(G)} 
MU_{*-|W|}\left( \prod_{V \in Irr^*(G)} 
BU(\nu_W(V)) \right).$$
Define the homomorphism $\vphi \colon \Omega^{U,G}_* \to F_*$
as sending a class $[M] \in \Omega^{U,G}_n$ to the class represented
by $M^G$ with reference map which 
classifies its normal bundle.
\end{definition}

This geometric picture of taking fixed sets of $G$-actions on manifolds
fits nicely with the homotopy theoretic picture we have been developing 
so far.

\begin{proposition}[tom Dieck]\label{P:compare}
The following diagram commutes
$$
   \begin{CD}
        \Omega^{U,G}_* @>\vphi>> F_* \\
        @VVPTV    @VViV\\
        MU^G_* @>\lambda>> (\Phi^G MU^G)_*,
   \end{CD}
$$
where $i$ is the inclusion map $$\bigoplus_{W \in R^+(G)} 
MU_{*-|W|}\left( \prod_{V \in Irr^*(G)} 
BU(\nu_W(V)) \right) \to \bigoplus_{W \in R(G)} MU_{*-|W|} 
\left( \prod_{Irr^*(G)} BU \right).$$
\end{proposition}

We may now compute the images of geometric classes in $MU^G_*$ under
localization by geometric means.

\begin{proposition}\label{P:Pninloc}
Let $V$ be an irreducible representation of $G$.
The image of $[\proj(n \oplus V)]$ in 
$(\Phi^G MU^G)_*$ is $Y_{n, V} + X$, where $X$ is 
$(e_{V^*})^{-n}$ for one-dimensional $V$ and is zero otherwise.
\end{proposition}

\begin{proof}
We use homogeneous coordinates.  There are two possible components
of the fixed sets.   The points whose coordinates ``in $V$'' are
zero, constitute a fixed $\cplx\proj^{n-1}$, whose normal bundle is
the tautological line bundle over $\cplx\proj^{n-1}$ tensored with $V$.  
As a class in $(\Phi^G MU^G)_*$, this manifold with reference map to 
$\prod_{Irr^*(G)} BU$ represents $Y_{n, V}$.
Alternately, when all other 
coordinates are zero the resulting submanifold is
the space of lines in $V$, which is an isolated fixed point when 
$V$ is one-dimensional and is a projective space with no fixed 
points, as $V$ has no non-trivial invariant subspaces, when $V$ has higher
dimension.  As classes in $(\Phi^G MU^G)_*$, isolated fixed points
represented negative powers of Euler classes.
\end{proof}

We finish this section by proving that for nice groups,
the inclusion of $R_0$ into $MU^G_*$
becomes an isomorophism after inverting Euler classes.

\begin{proof}[Proof of \refT{locisom}]
Recall that from \refC{comploc} we have that for nice groups 
$S^{-1}MU^G_* \cong MU_*[e_V, e_V^{-1}, Y_{i,V}].$ It suffices to consider
the image of $S^{-1}R_0$ in this ring.  The Euler classes and their inverses
are in this image by definition.  And by \refP{Pninloc}, the classes
$Y_{i,V}$ are in this image modulo negative powers of Euler classes.
\end{proof}

\section{Computation of $MU^G_*$}

By \refT{locisom}, for nice groups $G$ 
any $x \in MU^G_*$ can be multiplied by an
Euler class to get a class in $R_0$ modulo the annihilator of some Euler
class.  Our plan, which we carry out for abelian groups, 
is to build $MU^G_*$ from $R_0$ by division by 
Euler classes.  We are faced with two questions: ``when can one divide by 
an Euler class?'' and ``what are annhilators of Euler classes?''
When $G$ is a torus, \refT{seq} answers both of these questions.
Recall that $K(V)$ is the subgroup of $T$ which acts trivially on $V$. 

\begin{thmnonum}[Restatement of \refT{seq}]
The sequence
$$ 0 \to MU^T_* \overset{\cdot e_V}{\to} MU^T_* \overset{res^T_{K(V)}}{\to}
MU^{K(V)}_* \to 0 $$
is exact.
\end{thmnonum}

\begin{proof}
We construct the appropriate Gysin sequence and show that it breaks up
into short exact sequences.
 
Apply $\wt{MU_T}^*$ to the cofiber sequence
$ S(V)_+ \overset{i}{\to} S^0 \overset{j}{\to} S^V$
to get the long exact sequence
$$ \cdots \to \wt{MU_T}^{2n}(S^V) \overset{j^*}{\to}  {MU_T}^{2n}
\overset {i^*}{\to} MU_T^{2n}(S(V)) \overset{\delta}{\to}   
\wt{MU_T}^{2n+1}(S^V) \to \cdots. $$

As $MU_T$ has suspension isomorphisms for any representation, 
$\wt{MU_T}^{k}(S^V) \cong MU_T^{k - V}$.  
By \refP{period}, $MU_T^{k - V} \cong MU_G^{k-2}$.
The map $j^*$ is by definition multiplication by $e_V$.  

To compute $MU_T^k(S(V))$, we note that for a non-trivial  
irreducible representation of a torus $S(V)$ is homeomorphic to the
orbit space $T / K(V)$. 
But maps from this orbit space to $MU_T$ are
in one-to-one correspondence with maps from a single point in the
orbit to the $K(V)$-fixed set of $MU_T$, which is homeomorphic to 
the $K(V)$-fixed set of $MU_{K(V)}$.  We deduce that $MU_T^k(S(V))
\cong MU_{K(V)}^k$ and that $i^*$ is the restriction map.

By Comeza\~na's theorem (\refT{comez}), both $MU_T^*$ and $MU_{K(V)}^*$
are concentrated in even degrees.  Hence the long exact sequence above
yields the short exact of the theorem.
\end{proof}

\begin{remark}
Let $T = S^1$ and $\rho$ be the standard representation, so the restriction
to the kernel of $\rho$ is the augmentation map.
There is a geometric construction which reflects the fact that, by \refT{seq},
the augmentation ideal is principal, generated by $e_\rho$.

Let $f\colon X \to Y$ be an $S^1$-equivariant map of based spaces
which is null-homotopic upon forgetting the
$S^1$ action.  Let $F \colon X \times I \to Y$ be a null-homotopy.
Construct an $S^1$-equivariant map $f_{\Sigma(F)} \colon
X \times I \times S^1 \to Y$ by sending
$$ (x, t, \zeta) \mapsto \zeta \cdot F(\zeta^{-1} \cdot x, t).$$
This map passes to the quotient 
$$X \times I \times S^1 / \left( \{X \times 0 \times S^1\} \cup 
\{X \times 1 \times S^1\} \cup \{* \times I \times S^1\} \right),$$
which is $S^\rho \wedge X$.  When restricted to $S^0 \wedge X \subset
S^\rho \wedge X$ this map coincides with the orginal $f$, and thus
gives a ``quotient'' of $f$ by the class $S^0 \hookrightarrow S^\rho$.
\end{remark}

As in the introduction, fix a splitting $s_V$ of ${res}_{K(V)}$ as a 
map of $MU_*$-modules.  Let $\beta_V = s_V \circ {res}_{K(V)}$.
And for any $x \in MU^T_n$ let $\Gamma_V(x)$ be the unique class in
$MU^T_{n+2}$ such that $e_V \cdot \Gamma_V(x) = x - \beta_V(x)$.
The existence and uniqueness of $\Gamma_V(x)$ follow from \refT{seq}
and the fact that $x - \beta_V(x)$ is in the kernel of $res_{K(V)}$.

We are ready to prove our main theorem.

\begin{thmnonum}[Restatement of \refT{main}]
For any choice of splittings $s_V$, $MU^T_*$ is generated as an $MU_*$-algebra 
over the operations $\Gamma_V$ 
by the classes $e_V$ and $[\proj(\n \oplus V)]$, where $n \in \nats$ and
in all instances $V$ ranges over non-trivial
irreducible  complex representations of $T$.
Relations are as follows:
\begin{enumerate}
\item $e_V \Gamma_V(x) = x - \beta_V(x)$ \label{defgam}
\item $\Gamma_V(\beta_V(x)) = 0$ \label{gammabeta}
\item $\Gamma_V(e_V) = 1$ \label{gammaeul}
\item $\Gamma_V(xy) = \Gamma_V(x) y + \beta_V(x) \Gamma_V(y) + 
\Gamma_V(\beta_V(x) \beta_V(y))$  \label{product}
\item $ \Gamma_V \Gamma_W x  = \Gamma_W \Gamma_V x   + 
\Gamma_W \Gamma_V \beta_W(x)  - \Gamma_W \Gamma_V (e_W) \beta_V(\Gamma_W x) \\
- \Gamma_W \Gamma_V (\beta_V(e_W) \beta_V(\Gamma_W x)),$ \label{composition}
\end{enumerate} 
where $V$ and $W$  range over non-trivial irrreducible representations of
$T$ and $x$ and $y$ are any classes in $MU^G_*$.
\end{thmnonum}

\begin{proof}
By \refT{floc}, any class in $MU^T_*$ can be multiplied by an Euler class
to give a class in $R_0$ modulo the kernel of the canonical map from 
$MU^T_*$ to $S^{-1}MU^T_*$, where $S$ is the multiplicative set of non-trivial
Euler classes.  
By \refT{seq}, the kernel of the
map from $MU^T_*$ to the ring one obtains by inverting a single Euler class
is injective, so it follows that the map to $S^{-1}MU^T_*$ is injective.
Hence, any class in $MU^T_*$ is the quotient of some class in $R_0$ by
an Euler class.

Thus, we may filter $MU^T_*$ exhaustively as
$$R_* = R^0_* \subset R^1_* \subset \cdots \subset MU^T_*,$$
where $R^i_*$ is obtained from by adjoining to $R^{i-1}_*$ all $x \in MU^T_*$
such that $x \cdot e_V = y \in R^{i-1}_*$ for some $V$.  By \refT{seq} the
set of all such $y$ for a given $V$
is $\text{Ker}(res_{K(V)}) \cap R^{i-1}_*$.  The kernel of the
restriction map is clearly generated by all classes $y - \beta_V{y}$.
So we may obtain $R^i_*$ from $R^{i-1}_*$ by applying $\Gamma$ to every class
in $R^{i-1}_*$, which proves that $MU^T_*$ is generated over the operations
$\Gamma^V$ by $R_0$. 

Next, we note that the relations are readily verifyable.  
Relation \ref{defgam} holds by definition. And we may use the fact that 
multiplication by non-trivial
Euler classes is a monomorphism to verify relations \ref{gammabeta},
\ref{gammaeul} and  \ref{product} by
multiplying them by $e_V$, and \ref{composition} by multiplying it
by $e_V e_W$.

We are left to show that the relations are complete.  
%Our arguments will
%be straightforward induction using the filtration of $MU^T_*$ by the
%rings $R^i_*$.  
For convenience,
if $I = V_1, \cdots, V_k$ is a $k$-tuple of representations let 
$\Gamma_I(x) = \Gamma_{V_k} \Gamma_{V_{k-1}} \cdots \Gamma_{V_1} x$. 
We fix an ordering the representations of $T$.  

We claim that a multiplicative generating set for $MU^T_*$ is given by
classes $\Gamma_I(r)$ where $r$ is a generator of $R_0$
and $I$ is a (possibly empty) $k$-tuple of 
representations which respects the ordering we have imposed, as well
as classes $\Gamma_I(\prod \beta_{V_k}x_i)$, where $V_k$ is the 
minimimal representation in $I$.  
%Call this generating set $F$, and assume
%inductively that $F \cap R_i$ is a multiplicative generating set for 
%$R^{i-1}_*$.  
By the relation \ref{product}, to construct a generating set 
it suffices to consider classes $\Gamma_V(x)$, where 
$x$ is either primitive itself or a product of classes in the
image of $\beta_V$.  And by relation \ref{composition}, it suffices to 
consider within those classes only the ones $\Gamma_V(\Gamma_I(x))$ where
$V$ is greater than any of the representations in $I$.  

Next we give an additive basis for $MU^T_*$.  Fix an ordering on the 
generators of $R_0$ in which $r_i < r_j$ if $r_j$ is in the image of some
$\beta_{V}$ where $V$ is less than any $W$ such that $r_i$ is
in the image of $\beta_W$.  We define an additive basis for $MU^T_*$ in two
families.
Basis elements in the first family are 
the monomials $\Gamma_I(r_0) m$, where $m = \prod r_i$ is a monomial in $R_0$,
$I$ respects our ordering on 
the representations of $T$, 
$r_0$ is a generator of $R_0$ which is not in the image of $\beta_{V_k}$ where
$V_k$ is the minimal element of $I$, under our ordering $r_0$
is greater than any of 
the $r_i$, and for each representation $V \in I$ and each $i$ we
have that $e_V \neq r_i$.  Basis elements in the second family are 
$\Gamma_I(\prod r_i)$, where $I$ respects our ordering on the representations
of $T$ and each $r_i$ is a generator of $R_0$ in the
image of $\beta_{V_k}$ where $V_k$ is the minimal element of $I$.  

We check that this basis is linearly independent by mapping to $S^{-1}MU^T_*$.
Define a multiplicative basis of $S^{-1}MU^T_*$ using the images of elements
of $R_0$ along with the multiplicative inverses of Euler classes.  Extend our
ordering of generators of $R_0$ to an ordering of generators of 
$S^{-1}MU^T_*$ in which the inverses of Euler classes are less than any
generator of $R_0$.  Now order the monomials in $S^{-1}MU^T_*$ by a 
dictionary ordering.  Then the image of an additive basis element as
defined above 
$$\prod r_i \cdot \prod_{V_i \in I} e_{V_i}^{-1} 
+ {\text{lower order terms}}.$$
These images are clearly linearly independent.

Finally, we show that the relations suffice to reduce any product of 
multiplicative generators to a linear combination of additive basis
elements.  Let $m$ be a monomial in the multiplicative generators defined
above.  If $m = e_V \cdot \Gamma_I(x) \cdot \cdots,$ with $V \in I$, we
may use relation \ref{composition} to express $\Gamma_I(x)$ as a sum of
terms $\Gamma_V(\Gamma_J(x_i))$ and then use relation \ref{defgam} to
simplify.  We may thus reduce so that for each representation $V$ we
do not have both $e_V$ and $\Gamma_V$ appearing in $m$.  Next, note that
\ref{product} gives rise to the following relation
\begin{equation}\label{fromprod}
 x \Gamma_V(y) = \Gamma_V(x) (y - \beta_V(y)) + \beta_V(x) \Gamma_V(y).
\end{equation} 
We may use this relation repeatedly so that all of the operations
which appear in $m$ are applied to a single generator of our choosing.  
So we reduce to terms
of the form $\Gamma_I(r_0) \prod r_i$, where $r_0$ is greater than any $r_i$ 
in our ordering of generators of $R_0$.  Finally, we use the relation 
\ref{composition} to reorder the representations which appear in $I$.
Note that when we do so we may get terms $\Gamma_I (\prod \beta_V(x_i)) y$
which violate one of our conventions in that $\prod \beta_V(x_i)$ 
could be less
than some generators which appear in $y$.  We may then use equation 
\ref{fromprod} so that the operations are being applied to a maximal 
generator, followed by relation \ref{composition} to reorder.  This process
terminates.  At each stage we may associate a monomial in $R_0$ to a product
of our generators of $MU^T_*$ by forgetting all operations $\Gamma_V$.  
After an application of the relation \ref{composition} and equation 
\ref{fromprod} this associated monomial will be strictly smaller for
each term than the associated monomial for the original product.  
Once the associated monomials are of the minimal 
form $\prod \beta_{V_k}(x_i)$ 
where $V_k$ is minimal among representations appearing in the indexing
set for operations, we may use \ref{product} to equate the term with an
additive basis element in the second family.
\end{proof}

\section{The Completion Map and a Construction of Conner and Floyd}

From our computations, it is clear that the structure of $MU^G_*$ 
is governed by the operations $\Gamma_V$.  We call these operations
Conner-Floyd operations because in the special case of $G = S^1$,
$V = \rho$ the standard representation, and $[M]$ is a geometric class,
there is a construction of $\Gamma_\rho([M])$, which dates back to work
of Conner and Floyd.

\begin{definition}\label{D:gam}
Define $\gamma(M)$ for any stably complex $S^1$-manifold to be the
stably complex manifold
$$\gamma(M) = M \times_{S^1} S^3 \sqcup (-\overline{M}) \times \proj (1 \oplus
\rho),$$ where $S^3$ has the standard Hopf $S^1$-action,
$-\overline{M}$ is the $S^1$-manifold obtained from $M$ by imposing a trivial
action on $M$ and taking the opposite orientation, 
and the $S^1$-action on $M \times_{S^1} S^3$ is given by
\begin{equation}\label{E:dgam}
 \zeta \cdot \left[ m, z_1, z_2 \right] = \left[ \zeta \cdot
      m, z_1, \zeta z_2 \right].
\end{equation} 
Inductively define $\gamma^i(M)$ to be $\gamma(\gamma^{i-1}(M))$,
where $\gamma^0(M) = M$.
\end{definition}

\begin{proposition}\label{P:model}
Let $\rho$ be the standard representation of $S^1$.  And let $M$ be 
a stably complex $S^1$-manifold.  Then $\Gamma_\rho [M] = [\gamma(M)]$.
\end{proposition}

\begin{proof}
As the localization map is injective 
it suffices to check the equality in $S^{-1}MU^{S^1}_*$.  By \refP{compare} we
can compute the image of $[X]$, $[\overline{X}]$ and $[\Gamma(X)]$
in the localization at a full set of Euler classes by computing
fixed sets with normal bundle data.  The result follows easily as
the fixed sets of $\Gamma(X)$ are those of $X$ crossed with 
$\rho$ (when in the notation of equation \ref{E:dgam} above, $m$ is
fixed and $z_2 = 0$) 
along with an $\overline{X}$ crossed with $\rho$ (when $z_1 = 0$).  
In the localization, crossing with $\rho$ coincides with 
multiplying by $e_\rho^{-1}$.
\end{proof}

This geometric construction of a single Conner-Floyd operation allows us 
explicit understanding of the most important representation of $MU^T_*$,
namely the map from $MU^T_*$ to its completion at its augmentation ideal.
As a special case of \refT{seq}, we know that the augmentation ideal
of $MU^{S^1}_*$ is principal, generated by $e_\rho$.  Because the augmentation
map is split and multiplication
by $e_\rho$ is a monomorphism, the completion of $MU^{S^1}_*$ at its 
augmentation ideal is a power series ring over $MU_*$ where $e_\rho$ maps
to the power series variable under completion.  As an immediate
consequence of \refP{model} we have the following.

\begin{theorem}\label{T:gams}
Let $[M]$ be class in $MU^{S^1}_*$ which is the image under the 
Pontrjagin-Thom map of the class in geometric bordism represented by 
the complex $S^1$-manifold
$M$.  The image of $\left[ M \right]$ under the map from $MU^{S^1}_*$
to its completion at its augmentation ideal, which is isomorphic to
$MU_*[[x]]$, is the power series
$$ [\alpha(M)] + [\alpha(\gamma(M))]x + [\alpha(\gamma^2(M))]x^2
+ \cdots, $$
where $\alpha(\gamma^i(M))$ is the manifold obtained from 
$\gamma^i(M)$ simply by forgetting the $G$-action.
\end{theorem}

Understanding this completion map for geometric classes is important for
some geometric applications.  As mentioned in the introduction, 
Comeza/~na and May have proved that for abelian $G$, $(MU^G_*)_{\hat{I}}(X)
\cong MU_*(X \times_G EG)$.  So this completion homomorphism a the connection
between $MU^G$ and any cohomology theory which uses the Borel construction.
For example, let $\epsilon$ be a genus, that
is a ring homomorphism from $MU_*$ to some ring $E_*$.  For an 
extensive introduction to genera, see \cite{Seg1}.  
We may extend $\epsilon$ to an equivariant genus $\Omega^{U, G}_* \to
H^*(BG) \hat{\otimes} E_*$.  
Given a $G$-manifold $M$, take $M \times_G EG$
and use the genus $\epsilon$ to produce a class in $H^*(BG) \otimes E_*$
by ``integration over the fiber''.  
In our setting, for $G = S^1$, we may define this equivariant genus by taking
the image of a class $[M]$ under completion, namely 
$f \in (MU^{S^1}_*)_{\hat{I}} \cong MU^*[[x]]$, and applying $\epsilon$ 
term-wise to get a class in $H^*(\cplx \proj^{\infty}) \hat{\otimes} E_* \cong
E_*[[x]]$.

A genus $\epsilon$ 
is strongly multiplicative if for any fiber bundle of stably
complex manifolds $F \to E \to B$, 
$\epsilon(E) = \epsilon(F) \cdot \epsilon(B)$.  The following theorem
is a fundamental starting point in the study of genera, saying
essentially that strongly mutliplicative genera are rigid.

\begin{theorem}
Let $\epsilon$ be a strongly multiplicative genus.  Then for any class
$[M] \in \Omega^{U, S^1}_*$, the equivariant extension $\epsilon([M])$
is equal to $\epsilon([\alpha(M)]) \in E_* \subset E_*[[x]]$.
\end{theorem}

\begin{proof}
By \refT{gams} the image of $[M]$ under completion is 
$$ [\alpha(M)] + [\alpha(\gamma(M))]x + [\alpha(\gamma^2(M))]x^2
+ \cdots. $$
For any $X \in \Omega^{U, S^1}_*$ we have that 
$\epsilon([\alpha(\gamma(X))] = 0$ because $\epsilon$ is strongly
mulitplicative and by definition $\gamma(X)$ is the difference between
a twisted product and a trivial product of $X$ and $\cplx \proj^1$. 
\end{proof}

Returning to computation of the completion map on $MU^T_*$, we now
focus on Euler classes.

\begin{proposition}\label{P:euls}
The image of the Euler class $e_{\rho^{\otimes n}}$ in the completion
$(MU^{S^1}_*)_{\hat I}$
is $[n]_F x$, the $n$-series in the formal group law over $MU_*$.
\end{proposition}

\begin{proof}
As the map from $MU^G_*$ to its completion is a map of complex-oriented
equivariant cohomology theories, the Euler class of the bundle $V$ over a
point gets mapped to the Euler class of $V \times_G E_G$ over $B_G$.
For $G = S^1$, $V = \rho^{\otimes n}$ the resulting bundle is the $n$th-tensor
power of the tautological bundle over $B_{S^1}$, whose Euler class is
by definition the $n$-series.
\end{proof}

We are now ready to state our theorem about the image of the completion
map for $MU^T_*$.
When $T = (S^1)^k$, the completion of $MU^T_*$ at its augmentation ideal
is isomorphic to $MU_*[[x_1, \cdots x_k]]$.

\begin{definition}
Let $Y_n(x) \in MU_*[[x]]$ be the image of the class
$[\proj(\n \oplus \rho)]$ under the completion map.
\end{definition}

\begin{theorem}\label{T:oldmain}
Let $E$ be the set of all series 
$$[m_1]_F x_1 +_F \cdots +_F [m_k]_F x_k \in MU_*[[x_1, \cdots x_k]].$$
The image of
$MU^T_*$ in its completion at the augmentation ideal is contained in
the minimal
sub-ring $A$ of $MU_*[[x_1, \cdots, x_k]]$ which satisfies the
following two properties:
\begin{itemize}

\item $E \subset A$, and $A$ contains the series $Y_i(f)$ where $f \in E$
and $Y_i(f)$ are defined above.

\item If $\alpha f \in S$ then $\alpha \in A$, for any $f \in E$.

\end{itemize}

\end{theorem}

We can recover the image of $MU^G_*$ in its completion at the augmentation
ideal for general $G$ by reducing
$MU_*[[x_i]]$ modulo the ideal $([d_i]_F x_i)$, where $d_i$ are the orders
of the cyclic factors of $G$.

\begin{proof}
The first condition on $A$ says that the image
contains all images of classes in $R_0$.  Indeed, $E$ is the image of the
Euler classes.  And we check that the image of
$[\proj(i \oplus \rho^{\otimes n})]$ in $(MU^{S^1}_*)_{\hat{I}}$
is $Y_i([n]_F x)$, which follows from the fact that 
the $S^1$ action on $[\proj(i \oplus \rho^{\otimes n})]$
is pulled back from the $S^1$ action on $[\proj(i \oplus \rho)]$ by the
degree $n$ homomorphism from $S^1$ to itself.  By \refT{floc} we may build
any class in $MU^T_*$ by dividing classes in $R_0$ by Euler classes.
The second condition on $S$ accounts for all possible quotients by
Euler classes in the image.
\end{proof}

Suppose $f = a_0 + a_1 x + a_2 x^2 + \cdots$ 
is the image of $x \in MU^{S^1}_*$ under completion.  
Then the image of $\Gamma_\rho(x)$ under completion is
that $a_1 + a_2 x + a_3 x^2 + \cdots$ is in the image.
More generally, any $a_i + a_{i+1} x + \cdots$ is in the image of the
completion map.
So the property of a series being in the image of the completion map
depends only on the tail of the series.  It would be interesting to find
an ``analytic'' way to define this image.

\section{Applications and Further Remarks}

In this section we give an assortment of applications and indicate directions
for further inquiry.

Our first application is in answer to a question posed by Bott.  Suppose
a group acts on a manifold compatible with a stably complex structue
and that the fixed points of the action are isolated.  What can one say 
about the representations which appear as tangent spaces to the fixed points?
If there are only two fixed points, the representations must be dual, which
one can prove by Atiyah-Bott localization.  What happens for
three or more fixed points is an active area of inquiry in equivariant 
cohomology.  With our bordism techniques, we can get answer some of these
questions, as  wel as go beyond local information.

\begin{thmnonum}[Restatement of \refT{appl}]
Let $M$ be a stably-complex four dimensional
$S^1$-manifold with three isolated fixed points.  
Then $M$ is equivariantly cobordant to 
$\proj(\underline{1}\oplus V \oplus W)$ for some distinct non-trivial
irreducible representations $V$ and $W$ of $S^1$.
\end{thmnonum}

\begin{proof}
For convenience, let refer to the Euler class $e_{\rho^{\otimes n}} \in
MU^{S^1}_*$ by $e_n$.
A complex $S^1$ manifold $M$ with three isolated fixed points defines a 
class in $MU^{S^1}_*$ whose image under $\lambda \colon MU^{S^1}_* \to
S^{-1}MU^{S^1}_*$ is
$$\lambda([M]) = e_a^{-1} e_b^{-1} + e_c^{-1} e_d^{-1} + e_f^{-1}e_g^{-1}
$$ 
for some integers $a, \cdots, g$.  
We let $T$ denote 
$$e_a \cdots e_g \cdot \lambda[M] =
e_c e_d e_f e_g + e_a e_b e_f e_g + e_a e_b e_c e_d \in MU^{S^1}_*
\subset S^{-1}MU^{S^1}_*.$$  
%We say that a class in a localization is integral when it
%is in the image of the canonical map to that localization.

Without loss of generality, assume $a$ is greatest of the integers $a,
\cdots, g$
in absolute value.  As $T$ is divisible by $e_a$ in $MU^{S^1}_*$, 
\refT{seq} implies that $T$ restricted to 
$MU^{\ints/a}_*$ must be zero.  The Euler class $e_n$ restrict non-trivially 
to $MU^{\ints/a}_*$ unless $a|n$.  Therefore one of $c, d, f, g$, say $c$ must
be equal to $\pm a$.  We first claim that this number must be $-a$.  

Let $S_{\hat{a}}$ denote the multiplicative set generated by all
the Euler classes associated to irreducible representations except for $e_a$.
By localizing the modules in \refT{floc} and \refT{seq}, we find that
$S_{\hat{a}}^{-1}MU^T_*$ is generated over the operation $\Gamma_a$
by $S_{\hat{a}}^{-1} R_0$.
Suppose that $|b|, |d|, |f|, |g| < |a|$ and that $c = a$.  Then 
$$e_a^{-1}e_b^{-1} + e_a^{-1}e_d^{-1} + e_f^{-1}e_g^{-1}$$ is in the
image of the canonical map from $S_{\hat{a}}^{-1}MU^{S^1}_*$ to 
$S_{\hat{a}}^{-1} R_0$, as it is actually in the image of $\lambda$.  
Then  we must have that $e_b^{-1} + e_d^{-1}$ is divisible by $e_a$ and 
thus is zero in 
$S^{-1}MU^{\ints/a}_*$ where $S$ here is the multiplicative set of all
Euler classes of $\ints/a$.  This localization of $MU^{\ints/a}_*$ is
the the target of the restriction map 
from $S_{\hat{a}}^{-1}MU^{S^1}_*$.  And by abuse we are using the same
names for Euler classes for different groups. 
But because $|b|, |d| < |a|$, $e_b^{-1}$, $e_d^{-1}$ and their sum are
non-zero in $S^{-1}MU^{\ints/a}_*$.

It is straightforward to rule out cases where some of 
$|b|, |d|, |f|, |g|$ are equal to $|a|$.

Next, consider the class 
$$\lambda_{\hat{a}}([M]) - \lambda_{\hat{a}}([\proj(1 + \rho^{\otimes a})]) 
e_d^{-1}  \in S_{\hat{a}}^{-1}MU^{S^1}_*,$$ 
where $\lambda_{\hat{a}}$ is the canonical map to this localization.  
Its image under the map to the full localization is 
$$e_a^{-1}e_b^{-1} - e_a^{-1}e_d^{-1} + e_f^{-1}e_g^{-1},$$
which implies that $e_b^{-1} - e_d^{-1}$ is divisible by $e_a$ in 
$S_{\hat{a}}^{-1}MU^{S^1}_*$ or that $b \equiv d \pmod{a}$.
But because $|b|, |d| < |a|$ we have that $d = a \mp b$ depending on whether
$b$ is positive or negative. 

Finally, as $c = -a$ and $d = b-a$ consider $\lambda (
[M] - [\proj(1 \oplus \rho^{\otimes a} \oplus \rho^{\otimes b})])$, which 
will be equal $e_f^{-1}e_g^{-1} - e_{a-b}^{-1} e_{-b}^{-1}$.  Case analysis 
of necessary divisibilities as we have been doing 
implies that this difference must be zero, so that the fixed-set data of 
$[M]$ is isomorphic to that of 
$\proj(1 \oplus \rho^{\otimes a} \oplus \rho^{\otimes b})$.

Finally, by because the localization map $\lambda$ is injective,
this fixed-set data determines $[M]$ as 
in $S^1$-equivariant homotopical bordism uniquely, so that $[M]$ must equal 
$[\proj(1 \oplus \rho^{\otimes a} \oplus \rho^{\otimes b})]$ in $MU^{S^1}_4$.
But a theorem of Comeza\~na says that the Pontrijagin-Thom map from
$\Omega^{U,A}_*$ to $MU^A_*$ is injective for abelian groups $A$.  Hence
$M$ is cobordant to 
$\proj(1 \oplus \rho^{\otimes a} \oplus \rho^{\otimes b})$.

\end{proof}

Our next application answers a question about bordism of free 
$\ints/n$-manifolds posed to us by Milgram.  It is well-known that 
the spheres  $S(\oplus_k \rho^{\otimes m})$ for any $m$ relatively prime
to $n$ generate $MU_*(B_{\ints/n})$ as an $MU_*$-module.  How are these
bases related?

\begin{theorem}
Let $m$ and $n$ be relatively prime.
Let $Q(x)$ be a quotient of $x$ by $[m]_F x$ modulo $[n]_F x$ in $MU_*[[x]]$.
Define $a_i \in MU_*$ by $(Q(x))^k = a_0 + a_1 x + a_2 x^2 +\cdots$.
Then 
$$[S(\oplus_k \rho^{\otimes m})] = a_0 [S(\oplus_k \rho)] 
+ a_1 [S(\oplus_{k-1} \rho)] + \cdots + a_{k-1} [S(\rho)]$$
in $MU_*(B_{\ints/n})$.
\end{theorem}

\begin{proof}
We use an analog of the
simple fact that if $M$ is a $G$-manfold and $M \smallsetminus M^G$ has a
free $G$-action then $[\partial \nu(M^G)] = 0$ in $MU_*(B_G)$, where
$\partial \nu (M^G)$ is the boundary of a tubular neighborhood around
the fixed set $M^G$.  The null-bordism is defined by 
$M \smallsetminus \nu(M^G)$.  If
the fixed points of $M$ are isolated, this will give rise to a relation
among spheres with free $G$-actions.

Let $\alpha_0 = q^k$ where $q$ is 
a quotient of $e_\rho$ by $e_{\rho^{\otimes m}}$ in
$MU^{\ints/n}_*$.  Inductively, let $\alpha_i$ be a quotient of 
$\alpha_{i-1} - \overline{\alpha_{i-1}}$ by $e_\rho$ (note that this
quotient is not unique as we are working in $\ints/n$ equivariant
bordism.
Then the ``fixed sets'' of $\alpha_k$ are given by
$$\lambda(\alpha_k) = {e_{\rho^{\otimes m}}}^{-k} - \overline{\alpha_0}
{e_\rho}^{-k} - \overline{\alpha_1}{e_\rho}^{-k + 1} - \cdots -
\overline{\alpha_{k-1}}{e_\rho}^{-1}.$$
As ${e_V}^{-1}$ corresponds to a tubular neighborhood of an isolated fixed
point in geometric bordism, we can deduce via transversality arguments
for free $G$-actions that 
$$[S(\oplus_k \rho^{\otimes m})] - \overline{\alpha_0} [S(\oplus_k \rho)] 
- \overline{\alpha_1} [S(\oplus_{k-1} \rho)] - \cdots - 
\overline{\alpha_{k-1}} [S(\rho)] = 0$$ in $MU_*(B_{\ints/n})$.
But the image of $\alpha_0$ in $(MU^{\ints/n}_*)_{\hat I} \cong
MU_*[[x]]/[n]_F x$ is $(Q(x))^k$ from which we can read off that 
$\overline{\alpha_i} = a_i$.

Note that our expressions in $MU_*(B_{\ints/n})$ are independent of
the indeterminacy in choosing $q$ and the $\alpha_i$.
\end{proof}

This old idea of using $G$-manifolds to bound and thus give insight into
free $G$-manifolds has been codified by Greenlees's introduction of
local cohomology to equivariant stable homotopy theory \cite{Gr1}.  
%Essentially, $MU_*(B_G)$ is isomorphic to the
%cohomology of the Koszul-\v{C}ech complex
%$$\bigotimes_{V \in B} \left( MU^G_* \to MU^G_* \left[ \frac{1}{e_V} 
%\right] \right)$$
%where $S$ is a set of representations of $G$ such that $G$ acts freely
%on $\prod_{V \in B} S(V)$.  We have reproduced known $\text{Tor}$ classes 
%in $MU_*(B_{(\ints/2)^k})$ using this method.  
%Results along these lines
%will appear in \cite{Si2}.
Moreover, by work of Greenlees and May, 
the theories we have been studying provide
a unified framework in which to study the characteristic classes 
$E^*(BG)$ for any complex-oriented theory $E$.  We hope that our  
understanding of relevant commutative algebra can lead to new insights
into these techniques.

\end{document}